\documentstyle[12pt]{article}
\begin{document}

\title{A Recipe for Construction of the Critical Vertices for Left-Sector Stability of Interval Polynomials}
\author{Long Wang\thanks{Supported by National Key Project of China and Alexander von Humboldt Foundation of Germany. }\\Department of Mechanics and Engineering Science\\Peking University\\Beijing 100871\\China}
\date{}
\maketitle
\baselineskip 18pt
       
\begin{abstract}
For the left-sector stability of interval polynomials, it suffices to check a subset of its vertex polynomials. This paper provides a recipe for construction of these critical vertices. Illustrative examples are presented.

\end{abstract}

Left-sector stability is a very useful concept in characterizing the damping behavior of linear systems. Motivated by the seminal result of Kharitonov\cite{Khar}, a number of papers have concentrated on the left-sector stability of interval polynomials\cite{Soh,Shaw}. For the $n$-th order interval polynomials

\begin{equation}
{\bf K}(s)=\{f(s)\mid f(s)=a_0 + a_1 s + \ldots + a_n s^n,a_i\in
[a_i^- ,a_i^+ ],i=0,1,\ldots ,n\}
\end{equation}

\noindent and the left-sector stability region

\begin{equation}
  D = \left \{ \rho e^{j \theta} \mid \rho > 0\,,\,\frac{p \pi}{q} < \theta < 2 \pi - \frac{p \pi}{q} \right \}
\end{equation}

\noindent where $p\,,\,q$ are coprime positive integers and $\frac{1}{2} \leq \frac{p}{q} < 1\,,$ Soh and Foo\cite{Soh} have shown that all polynomials in ${\bf K}(s)$ have all their roots lying within $D$ if and only if $2q$ prespecified vertex polynomials in ${\bf K}(s)$ have all their roots lying within $D$. Some graphical methods\cite{Soh,Shaw} for determination of these critical vertices based on value set geometry\cite{Barm,Djaf,Fu} on the complex plane have also been proposed. In what follows, we will present a simple, direct procedure for construction of these critical vertices, which does not need graphical manipulation on the complex plane.

Any $f(s) \in {\bf K}(s)$ can be viewed as a vector

\begin{equation}
  f = [ a_0\,\, a_1\,\, \ldots\,\, a_n ]
\end{equation}

\noindent in the $(n+1)$-dimensional coefficient space. For any vertex polynomial

\begin{equation}
 f = [ a_0\,\, a_1\,\, \ldots\,\, a_n ]\,,\,a_i \in \{ a_i^- ,a_i^+ \}\,,\,i=0,1,\ldots ,n
\end{equation}

\noindent define the superscript of $f$ as

\begin{equation}
{\rm superscript}[f]=[{\rm superscript}[a_0]\,\,\, {\rm superscript}[a_1]\,\,\, \ldots\,\,\, {\rm superscript}[a_n]]
\end{equation}

\noindent For example, if

\begin{equation}
f=[ a_0^-\,\,\,a_1^+\,\,\,a_2^-\,\,\,a_3^-\,\,\,a_4^+\,\,\,a_5^+ ]
\end{equation}

\noindent then its superscript is

\begin{equation}
{\rm superscript}[f]=[ -\,\,\,+\,\,\,-\,\,\,-\,\,\,+\,\,\,+ ]
\end{equation}

\vspace{0.6cm}
\noindent {\sc Construction Procedure for Critical Vertices}:
\vspace{0.6cm}

\begin{description}

\item[Step 1.] Form a circle of $2q$ signs by placing positive sign ``$+$'' consecutively $q$ times and negative sign ``$-$'' consecutively $q$ times, respectively;

\item[Step 2.] Pick one of these $2q$ signs on the circle as the starting point;

\item[Step 3.] Move $n$ steps clockwisely on the circle with step length $p$ . Record the signs of the $n$ steps consecutively to form the superscript of a vertex polynomial;

\item[Step 4.] Change the starting point and repeat.

\end{description}

\vspace{0.6cm}
\noindent {\bf Theorem}

The $2q$ superscripts obtained by the above procedure are exactly the superscripts of the $2q$ critical vertices for the left-sector stability of interval polynomials.
\vspace{0.6cm}

\noindent Remark 1. The above construction procedure for the critical vertices does not need value set computation on the complex plane, thus is very suitable for computer implementation.

\noindent Remark 2. For low-order interval polynomials, the above procedure will automatically result in a reduced number of critical vertices.
\vspace{0.6cm}

\noindent Example 1.(Kharitonov's Theorem)

Let $\frac{p}{q} = \frac{1}{2}\,,$ which corresponds to the Hurwitz stability case. First, form a circle of signs with two positive signs and two negative signs consecutively; and then pick one sign on the circle as the starting point, move on the circle clockwisely with step length $1$, we get the superscript of one of the four Kharitonov polynomials

\begin{equation}
  [ +\,\,\,+\,\,\,-\,\,\,-\,\,\,+\,\,\,+\,\,\,-\,\,\,-\,\,\,\cdots \cdots ]
\end{equation}

\noindent By changing the starting point and repeating the procedure, we get the superscripts of the other three Kharitonov polynomials

\begin{equation}
  [ +\,\,\,-\,\,\,-\,\,\,+\,\,\,+\,\,\,-\,\,\,-\,\,\,+\,\,\,\cdots \cdots ]
\end{equation}

\begin{equation}
  [ -\,\,\,-\,\,\,+\,\,\,+\,\,\,-\,\,\,-\,\,\,+\,\,\,+\,\,\,\cdots \cdots ]
\end{equation}

\begin{equation}
  [ -\,\,\,+\,\,\,+\,\,\,-\,\,\,-\,\,\,+\,\,\,+\,\,\,-\,\,\,\cdots \cdots ]
\end{equation}
\vspace{0.6cm}

\noindent Example 2.

Let $\frac{p}{q} = \frac{3}{4}\,,$ form a circle of signs with four positive signs and four negative signs consecutively, pick one sign on the circle as the starting point, move on the circle clockwisely with step length $3$, we get the superscript of one of the eight critical vertex polynomials

\begin{equation}
  [ +\,\,\,+\,\,\,-\,\,\,+\,\,\,-\,\,\,-\,\,\,+\,\,\,-\,\,\,\cdots \cdots ]
\end{equation}

\noindent By changing the starting point and repeating the procedure, we get the superscripts of the other seven critical vertex polynomials

\begin{equation}
  [ +\,\,\,-\,\,\,-\,\,\,+\,\,\,-\,\,\,+\,\,\,+\,\,\,-\,\,\,\cdots \cdots ]
\end{equation}

\begin{equation}
  [ +\,\,\,-\,\,\,+\,\,\,+\,\,\,-\,\,\,+\,\,\,-\,\,\,-\,\,\,\cdots \cdots ]
\end{equation}

\begin{equation}
  [ +\,\,\,-\,\,\,+\,\,\,-\,\,\,-\,\,\,+\,\,\,-\,\,\,+\,\,\,\cdots \cdots ]
\end{equation}

\begin{equation}
  [ -\,\,\,-\,\,\,+\,\,\,-\,\,\,+\,\,\,+\,\,\,-\,\,\,+\,\,\,\cdots \cdots ]
\end{equation}

\begin{equation}
  [ -\,\,\,+\,\,\,+\,\,\,-\,\,\,+\,\,\,-\,\,\,-\,\,\,+\,\,\,\cdots \cdots ]
\end{equation}

\begin{equation}
  [ -\,\,\,+\,\,\,-\,\,\,-\,\,\,+\,\,\,-\,\,\,+\,\,\,+\,\,\,\cdots \cdots ]
\end{equation}

\begin{equation}
  [ -\,\,\,+\,\,\,-\,\,\,+\,\,\,+\,\,\,-\,\,\,+\,\,\,-\,\,\,\cdots \cdots ]
\end{equation}
\vspace{0.6cm}

\noindent Example 3.

Let $\frac{p}{q} = \frac{3}{5}\,,$ form a circle of signs with five positive signs and five negative signs consecutively, pick one sign on the circle as the starting point, move on the circle clockwisely with step length $3$, we get the superscript of one of the ten critical vertex polynomials

\begin{equation}
  [ +\,\,\,+\,\,\,-\,\,\,-\,\,\,+\,\,\,-\,\,\,-\,\,\,+\,\,\,+\,\,\,-\,\,\,\cdots \cdots ]
\end{equation}

\noindent By changing the starting point and repeating the procedure, we get the superscripts of the other nine critical vertex polynomials

\begin{equation}
  [ +\,\,\,+\,\,\,-\,\,\,+\,\,\,+\,\,\,-\,\,\,-\,\,\,+\,\,\,-\,\,\,-\,\,\,\cdots \cdots ]
\end{equation}

\begin{equation}
  [ +\,\,\,-\,\,\,-\,\,\,+\,\,\,+\,\,\,-\,\,\,+\,\,\,+\,\,\,-\,\,\,-\,\,\,\cdots \cdots ]
\end{equation}

\begin{equation}
  [ +\,\,\,-\,\,\,-\,\,\,+\,\,\,-\,\,\,-\,\,\,+\,\,\,+\,\,\,-\,\,\,+\,\,\,\cdots \cdots ]
\end{equation}

\begin{equation}
  [ +\,\,\,-\,\,\,+\,\,\,+\,\,\,-\,\,\,-\,\,\,+\,\,\,-\,\,\,-\,\,\,+\,\,\,\cdots \cdots ]
\end{equation}

\begin{equation}
  [ -\,\,\,-\,\,\,+\,\,\,+\,\,\,-\,\,\,+\,\,\,+\,\,\,-\,\,\,-\,\,\,+\,\,\,\cdots \cdots ]
\end{equation}

\begin{equation}
  [ -\,\,\,-\,\,\,+\,\,\,-\,\,\,-\,\,\,+\,\,\,+\,\,\,-\,\,\,+\,\,\,+\,\,\,\cdots \cdots ]
\end{equation}

\begin{equation}
  [ -\,\,\,+\,\,\,+\,\,\,-\,\,\,-\,\,\,+\,\,\,-\,\,\,-\,\,\,+\,\,\,+\,\,\,\cdots \cdots ]
\end{equation}

\begin{equation}
  [ -\,\,\,+\,\,\,+\,\,\,-\,\,\,+\,\,\,+\,\,\,-\,\,\,-\,\,\,+\,\,\,-\,\,\,\cdots \cdots ]
\end{equation}

\begin{equation}
  [ -\,\,\,+\,\,\,-\,\,\,-\,\,\,+\,\,\,+\,\,\,-\,\,\,+\,\,\,+\,\,\,-\,\,\,\cdots \cdots ]
\end{equation}

\end{document}